\documentclass[reqno, twoside,a4paper]{amsart}
\usepackage{style}
\renewcommand{\P}{\mathcal{P}}

\begin{document}

\title[Orthogonal polynomials]{
Entropy function and orthogonal polynomials
}
\author{R.\,V.\,Bessonov}

\address{
\begin{flushleft}
bessonov@pdmi.ras.ru\\\vspace{0.1cm}
St.\,Petersburg State University\\  
Universitetskaya nab. 7-9, 199034 St.\,Petersburg, RUSSIA\\
\vspace{0.1cm}
St.\,Petersburg Department of Steklov Mathematical Institute\\ Russian Academy of Sciences\\
Fontanka 27, 191023 St.Petersburg,  RUSSIA\\
\end{flushleft}
}

\thanks{
The work is supported by the Russian Science Foundation grant 19-71-30002.
}


\begin{abstract} 
We give a simple proof of a classical theorem by A.\,M\'at\'e, P.\,Nevai, and V.\,Totik on asymptotic behavior of orthogonal polynomials on the unit circle. It is based on a new real-variable approach involving an entropy estimate for the orthogonality measure. Our second result is an extension of a theorem by G.~Freud on averaged convergence of Fourier series. We also  discuss some related open problems in the theory of orthogonal polynomials on the unit circle. 
\end{abstract}

\maketitle

\section{Introduction}
Consider a probability measure $\mu = w\,dm+\mu_s$ on the unit circle $\T = \{z: |z| = 1\}$ of the complex plane, $\C$. Here $w \in L^1(\T)$  and the measure $\mu_s$ is singular with respect to the Lebesgue measure $m$ on $\T$. As usual, we use normalization $m(\T) = 1$ and write $L^1(\T)$ for $L^1(m)$. Assume that $\log w \in L^1(\T)$, and associate with $\mu$ the so-called Szeg\H{o} function,
$$
D_\mu: z \mapsto \exp\left(\frac{1}{2}\int_{\T}\log w(\xi) \frac{1+\bar\xi z}{1 - \bar \xi z}\,dm(\xi)\right), \qquad |z| <1,
$$
that is, the outer function in the open unit disk with $|D_\mu|^2 = w$ almost everywhere on $\T$ in the sense of non-tangential boundary values. Let $\phi_n$, $n \ge 0$, denote orthonormal polynomials in $L^2(\mu)$ obtained from harmonics $z^n$ via the Gram-Schmidt orthogonalization process. Define also the reflected polynomials $\phi_n^*: z \mapsto z^n\ov{\phi_n(1/\bar z)}$ associated with $\mu$, and note that $|\phi_n^*| = |\phi_n|$ on $\T$, $\phi_n^*(0) > 0$.  

\medskip

Given a measure $\lambda$ on $\C$, we say that  $\xi_0 \in \supp \lambda$ is a Lebesgue point of $f \in L^p(\lambda)$, $1 \le p < \infty$, if there exists a representative of $f$ (to be denoted by the same letter) such that
$$
\lim_{\eps \to 0}\frac{1}{\lambda(I_{\xi_0, \eps})}\int_{I_{\xi_0, \eps}}|f - f(\xi_0)|^p\,d\lambda = 0, \qquad I_{\xi_0, \eps} =\{\xi \in \C: \; |\xi - \xi_0| < \eps\}.
$$
Note that a point $\xi_0$ may or may not be a Lebesgue point of a function $f \in L^p(\lambda)$ depending on which value $p$ we consider. 

\medskip

Our aim here is to give a new proof of the following remarkable result due to A.~M\'at\'e, P.~Nevai, and V.~Totik \cite{MNT91}. 
\begin{Thm}\label{t1}
Let $\mu = w\,dm + \mu_s$ be a probability measure on $\T$. Suppose that $\xi_0 \in \T$ is a Lebesgue point of $w$, $\log w \in L^1(\T)$, and  $\lim_{\eps \to 0}\mu_s(I_{\xi_0, \eps})/m(I_{\xi_0, \eps}) = 0$. Then 
\begin{equation}\label{eq20}
	\lim_{n \to+\infty}\frac{1}{n}\sum_{k=0}^{n-1}|\phi_n(\xi_0)|^2 = w(\xi_0)^{-1}.
\end{equation}
\end{Thm}
Compared with the original statement from 1991, assumptions of Theorem \ref{t1} do not require $\xi_0$ to be a Lebesgue point of the Szeg\H{o} function $D_\mu \in L^1(\T)$. In fact, this extra assumption was proved to be superfluous in a recent work by V.\,Totik \cite{Tot16} (see Theorem 1.3 therein). 

\medskip

Our proof is essentially real-variable: the only tool from complex analysis we need is the Poisson formula for real-valued harmonic functions on $\T$. This makes our arguments close to S.~Khurshchev's proof \cite{KH01} of the classical Szeg\H{o} theorem. We also give an estimate for the rate of convergence in \eqref{eq20},  see Proposition~\ref{p1}.

\medskip

It worth be mentioned that even in the situation where $\mu$ has the form $\mu = w\,dm$ with a positive continuous weight $w$ on $\T$, we may have $\limsup|\phi_n(\xi_0)| = +\infty$ at some points $\xi_0 \in \T$ (in particular, one cannot omit Ces\`aro averaging in \eqref{eq20} in general). This result was proved by M.~Ambroladze \cite{Amb91} using ideas of E.~Rahmanov's work \cite{Rahmanov}, see also A.~Aptekarev,  S.~Denisov, D.~Tulyakov  \cite{ADT} for recent advances in estimating $\max_{\xi \in \T}|\phi_n(\xi)|$ as $n \to +\infty$.

\medskip

Historically, the interest in studying the average asymptotic behavior of orthogonal polynomials on the unit circle was strongly motivated by a classical result of G.~Freud \cite{Fr52} from 1952 (subsequent expositions can be found in Chapter IV of \cite{Freudbook} or in Section 4 of \cite{Nev86}). It can be formulated as follows. Let $\nu$ be a finite measure on the segment $[-1, 1]$. Consider orthonormal polynomials $p_n$, $\deg p_n = n$, $n \ge 0$, with respect to $\nu$. Given a function $f \in L^2(\nu)$, denote by $S_k(f, \cdot)$ the partial sums of its Fourier series with respect to the basis $\{p_n\}$:
$$
S_k(f, x) = \sum_{j = 0}^{k}c_j p_j(x), \qquad c_j = (f, p_j)_{L^2(\nu)}, \qquad x \in [-1,1].
$$
Suppose that $x_0 \in [-1,1]$ is such that 
\begin{equation}\label{eq18}
\sup\frac{1}{n}\sum_{k=0}^{n-1}|p_k(x_0)|^2 < \infty. 
\end{equation}
With this assumptions, Freud's theorem says: {\it if $x_0$ is a Lebesgue point of $f \in L^2(\nu)$, then the Fourier series of $f$ is strongly Ces\`aro summable at $x_0$ to $f(x_0)$, that is, }
$$
\lim_{n \to +\infty}\frac{1}{n}\sum_{k=0}^{n-1}|S_k(f, x_0) - f(x_0)| = 0.
$$
A combination of Freud's and M\'at\'e-Nevai-Totik's theorems implies that if $\nu = v\,dx+\nu_s$ is such that 
\begin{equation}\label{eq19}
\int_{-1}^{1}\frac{\log v(x)}{\sqrt{1 - x^2}}\,dx > -\infty,
\end{equation}
then each $L^2(\nu)$-convergent Fourier series $\sum_{k \ge 0}c_k p_k$  is strongly Ces\`aro summable Lebesgue almost everywhere on $[-1,1]$. A proof of this corollary can be found in A.\,Mate, P.\,Nevai \cite{MN80}.

\medskip

It is natural to expect that ideas of Freud's work \cite{Fr52} are applicable in the setting of orthogonal polynomials on the unit circle. As a result, we potentially should have the following statement: if $\mu = w\,dm + \mu_s$ is a probability measure on $\T$ with $\log w \in L^1(\T)$, then every $L^2(\mu)$-convergent Fourier series $\sum_{k \ge 0}c_k \phi_k$ is strongly Ces\`aro summable Lebesgue almost everywhere on $\T$. Surprisingly, an attempt to ``transfer'' Freud's proof from $[-1,1]$ to $\T$ encounters major difficulties, see M\'at\'e \cite{Mate}. In particular, aforementioned statement is not proved yet. The main difference between theories on $[-1,1]$ and $\T$ in this case is hidden in the fact that Christoffel-Darboux kernel for polynomials on $\T$ contains reflected polynomials $\phi_n^*$ that are far from being orthonormal system (for instance, $\phi_n^* = 1$ for all $n$ if $\mu = m$).

\medskip

As we will see in a moment, from the point of view of summability of Fourier series, a ``correct'' replacement of polynomials $\{p_n\}_{n\ge 0}$ in $L^2(\nu)$ is the so-called CMV basis $\{\chi_n\}_{n \ge 0}$ in $L^2(\mu)$. To be precise, fix a probability measure $\mu$ on the unit circle $\T$ supported on an infinite set (the latter implies that $L^2(\mu)$ is infinite dimensional and the convergence issues make sense). Define $\chi_n$, $n \ge 0$, to be the orthonormal sequence obtained by Gram-Schmidt orthogonalization in $L^2(\mu)$ of har\-mo\-nics $1, z, \bar z, z^2, \bar z^2, \ldots$ (in this order). In fact, we have
$$
\chi_{2k} = \bar z^{k}\phi_{2k}^{*}, \qquad  
\chi_{2k+1} = \bar z^{k}\phi_{2k+1}, \qquad z \in \T, \qquad
k \ge 0.
$$
These trigonometric polynomials were introduced into a wide use of orthogonal polynomial community by M.\,Cantero, L.\,Moral, and L.\,Vel\'azquez \cite{CMV}. The fact that $\{\chi_n\}_{n \ge 0}$ forms an orthonormal basis in $L^2(\mu)$ can be easily verified by using orthogonality of usual polynomials $\{\phi_n\}$ generated by $\mu$. Below we use notation
$$
\bS_k(f, \xi) = \sum_{j = 0}^{k}c_j \chi_j(\xi), \qquad c_j = (f, \chi_j)_{L^2(\mu)}, \qquad \xi \in \T,
$$
for the corresponding partial sums of the Fourier series of a function $f \in L^2(\mu)$. We also will need a notation for the Poisson extension of $\mu$,
$$
\P(\mu, z) = \int_{\T}\frac{1-|z|^2}{|1 - \ov{\xi} z|^2}\,d\mu(\xi), \qquad |z|<1.
$$
Our version of Freud's theorem for CMV basis $\{\chi_n\}_{n \ge 0}$ reads as follows.
\begin{Thm}\label{t2}
Let $\mu$ be a probability measure on $\T$. Suppose that $\xi_0 \in \T$ is such that
\begin{equation}\label{eq17}
\frac{1}{n}\sum_{k=0}^{n}|\chi_n(\xi_0)|^2 \le \frac{c}{\P(\mu, z_n)}, \qquad z_n = \left(1 - \tfrac{1}{n}\right)\xi_0,
\end{equation}
for all $n \ge 1$ and a constant $c$. Assume, moreover, that $\xi_0$ is a Lebesgue point of $f \in L^2(\mu)$. Then 
$$
\lim_{n \to +\infty}\frac{1}{n}\sum_{k=0}^{n-1}|\bS_k(f, \xi_0) - f(\xi_0)| = 0.
$$
In particular, if $\mu = w\,dm + \mu_s$ is such that $\log w \in L^1(\T)$, then the Fourier series of $f$ is strongly Ces\`aro summable Lebesgue almost everywhere on $\T$.
\end{Thm}
We end this section with two problems that seem to be important for later development of summability theory of general Fourier series.
\begin{Prob}\label{prob1}
Describe the class of probability measures on $\T$ such that \eqref{eq17} holds $\mu$-almost everywhere on $\T$.
\end{Prob}
To the author's knowledge, at the present moment there is no example of an infinitely supported measure $\mu$ on $\T$ that does not satisfy \eqref{eq17} $\mu$-almost everywhere. Moreover,  every measure $\mu$ of the form $\mu =  w\,dm + \mu_s$ with $w = 0$ or $\log w \in L^1(\T)$ and a purely discrete singular part $\mu_s$ satisfies \eqref{eq17} $\mu$-almost everywhere, see Proposition \ref{p4}. Thus, the class of measures in Problem \ref{t2} contains both very regular and very singular measures, as well as their linear combinations. To emphasize the depth of the problem, let us note that it took more than 25 years after G.\,Freud's work \cite{Fr52} until A.\,Mate and P.\,Nevai  \cite{MN80} proved that each measure satisfying Szeg\H{o} condition \eqref{eq19} also satisfies Freud's bound \eqref{eq18} Lebesgue almost everywhere on $[-1,1]$.
\begin{Prob}\label{prob2}
Let $\mu = w\,dm + \mu_s$ be a probability measure on $\T$, $\log w \in L^1(\T)$, and let $\{c_n\} \subset \C$ be such that $\sum_{n}|c_n|^2 < \infty$. Prove that $\sum_{n \ge 0}c_n \phi_n$ converges Lebesgue almost everywhere on $\T$ in Ces\`aro sense or construct a counterexample.
\end{Prob}
A particular unsolved case of Problem~\ref{prob2} concerns the situation when $c_n = \ov{\phi_n(0)}$. Then $\sum c_n\phi_n = D_{\mu}^{-1}\ov{D_{\mu}^{-1}(0)}$ in $L^2(\mu)$, and convergence results for this series have interesting applications in scattering theory, see Section \ref{s3}.

\section{Entropy function of a measure. Proof of Theorem \ref{t1}}\label{s2}

Let $\mu = w\,dm + \mu_s$ be a probability measure on $\T$. The Schur function of $\mu$ is the analytic function $f$ in the open unit disk, $\D = \{z \in \C: |z| < 1\}$, defined by the relation
\begin{equation}\label{eq00}
\frac{1 + zf(z)}{1-zf(z)} = \int_{\T}\frac{1 + \bar \xi z}{1 - \bar \xi z}\,d\mu(\xi), \qquad z \in \D.
\end{equation}
Taking the real parts, 
$$
\frac{1 - |zf(z)|^2}{|1-zf(z)|^2} = \int_{\T}\frac{1 - |z|^2}{|1 - \bar \xi z|^2}\,d\mu(\xi), \qquad z \in \D,
$$
we see that the function $f$ is contractive in $\D$. Set $f_0 = f$ and define the Schur iterates, $f_n$, of $f$: 
\begin{equation}\label{eq03}
zf_{n+1}(z) = \frac{f_n(z) - f_n(0)}{1 - \ov{f_n(0)}f_n(z)}, \qquad z \in \D, \qquad n \ge 0. 
\end{equation}
It can be seen that $|f_{n}| < 1$ in $\D$ and $f_{n+1}$ is correctly defined for all $n \ge 0$ unless $f$ is a finite Blaschke product (that is, $\mu$ is supported on a finite subset of $\T$ -- the case that we exclude from our consideration).
A measure $\mu = w\,dm+\mu_s$ is said to belong to the Szeg\H{o} class $\szc$ if $\log w \in L^1(\T)$. The Szeg\H{o} formula says
\begin{equation}\label{eq3}
\int_{\T}\log w\, dm = \log\prod_{k \ge 0}^{\infty}(1-|f_k(0)|^2).
\end{equation}
In particular, we have $\mu \in \szc$ if and only if $\sum_{k \ge 0}|f_k(0)|^2 < \infty$. Formula \eqref{eq3} was generalized in \cite{BD2020} in the following way. For $z \in \D$, denote 
\begin{align*}
\K(\mu, z) 
&= \log \int_{\T}\frac{1 - |z|^2}{|1 - \bar \xi z|^2}\,d\mu(\xi) - \int_{\T}\log w(\xi)\frac{1 - |z|^2}{|1 - \bar \xi z|^2}\,dm(\xi),\\
&= \log \P(\mu, z) - \P(\log w, z).
\end{align*}
Jensen's inequality implies that $\K(\mu, z) \ge 0$ for every $\mu \in \szc$, $z \in \D$, and $\K(\mu, z)= 0$ if and only if $\mu = m$. In fact, this function measures a  ``size of irregularity'' of $\mu$ near the point $z/|z|$: if $\K(\mu, z)$ is small, then $\mu$ is close to $\P(\mu, z)\, dm$ near $z/|z|$ (we assume that $z$ is close to the unit circle). This explains the name ``entropy function of $\mu$'' for $\K(\mu, z)$. Szeg\H{o}'s formula \eqref{eq3} can be rewritten in the form
$$
\K(\mu, 0) = \log\prod_{k \ge 0}^{\infty}\frac{1}{1-|f_k(0)|^2}.
$$ 
It was proved in \cite{BD2020} that
\begin{equation}\label{eq6}
\K(\mu, z) = \log\prod_{k \ge 0}^{\infty}\frac{1 - |z f_k(z)|^2}{1- |f_k(z)|^2}, \qquad z \in \D.
\end{equation}
Extension \eqref{eq6} of the Szeg\H{o} formula implies the following important fact. 
\begin{Lem}\label{l1}
Let $\mu$ be a measure in $\szc$, and let $f$ be its Schur function. Then
$$
(1 - |z|^2)\sum_{k \ge 0}^{\infty}\frac{|f_k(z)|^2}{1 - |f_k(z)|^2} \le e^{\K(\mu, z)} - 1, \qquad z \in \D.
$$	
\end{Lem}
\beginpf For every $z \in \D$ we have  
$$
\frac{1 - |z f_k(z)|^2}{1 - |f_k(z)|^2} = 1 + (1 - |z|^2)\frac{|f_k(z)|^2}{1 - |f_k(z)|^2}.
$$
Therefore, for every $n \ge 0$ we have
$$
1 + (1 - |z|^2)\sum_{k = 0}^{n}\frac{|f_k(z)|^2}{1 - |f_k(z)|^2} \le \prod_{k=0}^{n}\left(1 + (1 - |z|^2)\frac{|f_k(z)|^2}{1 - |f_k(z)|^2}\right) \le e^{\K(\mu, z)}. 
$$
Since $n$ is arbitrary, the lemma follows. \qed

\medskip

\begin{Lem}\label{l2}
For every $\mu \in \szc$ we have
$$
\frac{1}{n}\sum_{k = 0}^{n-1}|D_\mu(z)\phi_k^*(z)|^2 \le 1 + 2\sqrt{\frac{e^{\K(\mu, z)} - 1}{n(1-|z|^2)}} +
4\frac{e^{\K(\mu, z)} - 1}{n(1-|z|^2)}, \qquad z \in \D.
$$
\end{Lem}
\beginpf For every $z \in \D$ we have
$$
|\phi_n^*(z)D_{\mu}(z)|^2 \le \left| \int_{\T}\phi_n^*(\xi)D_{\mu}(\xi)\frac{1-|z|^2}{|1 - \bar\xi z|^2}\,dm(\xi) \right|^2 \le \int_{\T}|\phi_n^*(\xi)|^2\frac{1-|z|^2}{|1 - \bar\xi z|^2}\,d\mu(\xi),
$$
by Jensen's inequality and the fact that 
$
|D_\mu|^2\,dm = w\,dm \le w\,dm + \mu_s = \mu.
$ 
Set $b_n = \phi_n/\phi_n^*$. Khrushchev's formula (see Theorem 3 in \cite{KH01}) gives  
$$
\int_{\T}|\phi_n^*(\xi)|^2\frac{1-|z|^2}{|1 - \bar\xi z|^2}\,d\mu(\xi) = \frac{1 - |z b_n(z) f_n(z)|^2}{|1 - z b_n(z)f_n(z)|^2}.
$$
Since $b_n$ is a finite Blaschke product (it is well-known that polynomials $\phi_n^*$ do not have zeroes in $\D$, see, e.g., Section 1.7 in \cite{Simonbook}), we  have
$$
|\phi_n^*(z)D_{\mu}(z)|^2 \le \frac{1 + |f_n(z)|}{1 - |f_n(z)|} \le 1 +\frac{2|f_n(z)|}{\sqrt{1 - |f_n(z)|^2}} + \frac{4|f_n(z)|^2}{1 - |f_n(z)|^2}, 
$$ 
where we used the elementary inequality
$$
\frac{1+x}{1 - x} \le 1+\frac{2x}{\sqrt{1 - x^2}} + \frac{4x^2}{1 - x^2}, \qquad 0 \le x < 1.
$$
From here and the inequality $\frac{1}{n}\sum_{0}^{n-1} t_k \le \sqrt{\frac{1}{n}\sum_{0}^{n-1} t_k^2}$ we obtain
$$
\frac{1}{n}\sum_{k = 0}^{n-1} |\phi_k^*(z)D_{\mu}(z)|^2 \le 1 + 2\sqrt{\eta_n} + 4\eta_n, \qquad \eta_n = \frac{1}{n}\sum_{k = 0}^{n-1} \frac{|f_k(z)|^2}{1 - |f_k(z)|^2}. 
$$
It remains to use Lemma \ref{l1}. \qed

\medskip

The following lemma is is due to A.\,M\'at\'e and P.\,Nevai \cite{MN80}.

\begin{Lem}\label{l3}
If $p_n$ is a polynomial of degree at most $n$ without zeroes in $\D$, then 
$$|p(r\xi)| \ge \left(\frac{1+r}{2}\right)^n|p(\xi)|$$ 
for every $\xi \in \T$ and $0 \le r \le 1$. 
\end{Lem}
\beginpf Since the inequality we want to prove is multiplicative, one can assume that $p = \lambda - z$ for some $\lambda$ such that $|\lambda| \ge 1$. Then the expression
$$
\frac{|p(r\xi)|}{|p(\xi)|} = \frac{|\lambda - r\xi|}{|\lambda - \xi|}  
$$
attains its minimum at $\lambda = -\xi$, which proves the statement. \qed

\medskip

\noindent{\bf Proof of Theorem \ref{t1}.} Fix some $\delta \in (0, 1]$ and denote $r_n = 1-\delta/n$. We have
$n(1 - r_n^2) \ge n(1 - r_n) = \delta$ for all $n$. Choosing $z = z_n = r_n \xi_0$ in Lemma \ref{l2}, we get 
\begin{equation}\notag 
\frac{1}{n}\sum_{k=0}^{n-1}|D_\mu(z_n)\phi^*_k(z_n)|^2 \le 1 + 2\sqrt{\eps_n} + 4\eps_n, \qquad \eps_n = \frac{e^{\K(\mu,z_n)}-1}{\delta}.
\end{equation}
Recall that polynomials $\phi_n^*$ do not have zeroes in $\D$. Taking into account Lemma \ref{l3} and the definition of $D_\mu$, we can proceed as follows:
\begin{equation}\label{eq22}
\frac{1}{n}\sum_{k=1}^{n-1}e^{\P(\log w, z_n)} |\phi_k(\xi_0)|^2 \le \left(1 - \frac{\delta}{2n}\right)^{-2n} (1 + 2\sqrt{\eps_n} + 4\eps_n).
\end{equation}
Assumptions imposed on $\xi_0$ imply
$\P(\log w, z_n) \to \log w(\xi_0)$, $\P(w, z_n) \to w(\xi_0)$, $\K(\mu, z_n) \to 0$, and $\eps_n \to 0$ as $n \to +\infty$. Hence, 
$$
\limsup_{n \to +\infty} \frac{1}{n}\sum_{k=0}^{n-1}|\phi_k(\xi_0)|^2 \le e^{\delta} w(\xi_0)^{-1}.
$$
Since $\delta \in (0,1)$ is arbitrary, this completes the proof of the inequality ``$\le$'' in \eqref{eq20}. For the proof of the opposite inequality, we follow \cite{MN80}. Denote by $k_{P_{n-1}, \mu, \xi}$ the reproducing kernel at a point $\xi \in \C$ in the Hilbert space $P_{n-1}$ of polynomials of degree at most $n-1$ with the inner product inherited from $L^2(\mu)$. Then 
$$
\|k_{P_{n-1}, \mu, \xi}\|^{2}_{L^2(\mu)} = \sup\{|p(\xi)|^2: \; p \in P_{n-1}, \; \|p\|_{L^2(\mu)} \le 1\}.
$$
In particular, we have
$$
\|k_{P_{n-1}, \mu, \xi_0}\|^{2}_{L^2(\mu)} \ge \frac{|p_{n-1}(\xi_0)|^2}{\|p_{n-1}\|^{2}_{L^2(\mu)}}, \qquad p_{n-1}: z \mapsto \sum_{k = 0}^{n-1}\bar\xi_0^k z^k.
$$
This can be rewritten in the form
\begin{equation}\label{eq24}
\frac{\|k_{P_{n-1}, \mu, \xi_0}\|^{2}_{L^2(\mu)}}{n} \ge \frac{1}{\F_{n-1}(\mu, \xi_0)},
\end{equation}
where 
$$
\F_{n-1}(\mu, \xi_0) = \frac{1}{n}\int_{\T}|p_{n-1}(\xi)|^2\,d\mu(\xi)
$$
is the classical Fejer mean of order $n-1$ centered at $\xi_0$. Since $\xi_0$ is the Lebesgue point of $w$, we have $\lim_{n \to +\infty}\F_{n}(\mu, \xi_0) = w(\xi_0)$. 
On the other hand, $\{\phi_k\}_{k = 0}^{n-1}$ forms the orthonormal basis in $P_{n-1}$, hence
\begin{equation}\label{eq25}
\frac{\|k_{P_{n-1},\mu,  \xi_0}\|^{2}_{L^2(\mu)}}{n} = \frac{1}{n}\sum_{k=0}^{n-1}|\phi_k(\xi_0)|^2, 
\end{equation}
and the result follows. \qed

\medskip

It is well known that M\'at\'e-Nevai-Totik's theorem is strongly related to universality limits in the random matrix theory. On the level of orthogonal polynomials, one can say that universality holds at a point $x \in \R$ for the orthogonality measure $\nu = v\,dx+ \nu_s$ compactly supported on $\R$ if there exists the limit
\begin{equation}\label{eq27}
\frac{K_n\left(x+\frac{a}{v(x)K_n(x,x)}, x+\frac{b}{v(x)K_n(x,x)}\right)}{K_n(x,x)} \to  \frac{\sin\pi(a-b)}{\pi(a-b)}, \qquad n \to +\infty,
\end{equation}
uniformly in $a,b$ lying in some compact subset of $\R$. Here 
$$K_n(x,y) = \frac{1}{n}\sum_{0}^{n-1}p_k(x)p_k(y), \qquad x, y \in \R,$$ 
denotes the Christoffel-Darboux kernel generated by orthogonal polynomials $\{p_k\}$ of $\nu$. After D.\,Lubinsky's seminal paper \cite{Lub09}, there was a substantial progress in relaxing local conditions on $\nu$ under which universality holds. We mention here the works by E.\,Findley \cite{Fin08} and V.\,Totik \cite{Tot16} where local versions of the M\'at\'e-Nevai-Totik theorem were established with this aim. However, to the author's knowledge, no estimate of the rate of convergence in \eqref{eq20} and \eqref{eq27} is available at the present moment. Below we give such an estimate for \eqref{eq20} in terms of $\K(\mu,z)$. Recall that $\K(\mu,z)$ essentially depends on the local regularity properties of $\mu$. In particular, it seems plausible that the estimate in Proposition \ref{p1} below could be ``localized'' and used later to give an estimate for the rate of convergence in \eqref{eq27} in the general setting of universality considered in \cite{Tot16}.

\medskip

\noindent For $n\ge 1$, set 
\begin{align*}
\K_n(\xi_0) 
&= \sup_{\delta \in (0,1)}\K(\mu, (1- \delta/n)\xi_0), \\
\P_n(\xi_0) 
&= \inf_{\delta \in (0,1)}\P(\mu, (1- \delta/n)\xi_0), \\
\F_n(\xi_0) 
&= \F_{n-1}(\mu, \xi_0).
\end{align*}
Note that these quantities are completely determined by $\mu$, $n$, and $\xi_0$. We have $\K_n(\xi_0) \to 0$, $\P_n(\xi_0) \to w(\xi_0)$, and $\F_n(\xi_0) \to w(\xi_0)$ as $n \to +\infty$ at every Lebesgue point $\xi_0 \in \T$ of $w, \log w \in L^1(\T)$ such that $\lim_{\eps \to 0} \mu_s(I_{s,\eps})/\eps = 0$. 
\begin{Prop}\label{p1} Let $\mu$, $\xi_0$ be as in Theorem \ref{t1}. Then 
\begin{equation}\label{eq26}
\frac{1}{\F_n(\xi_0) } \le \frac{1}{n}\sum_{k = 0}^{n-1}|\phi_k(\xi_0)|^2 \le \frac{1}{\P_n(\xi_0)}
+64\frac{\sqrt[4]{\K_n(\xi_0)}}{\P_n(\xi_0)}
\end{equation}
for all $n \ge 1$ such that $\K_n(\xi_0) \le 1$. In other words, for such $n$ we have
\begin{equation}\notag
\left|\frac{1}{n}\sum_{k = 0}^{n-1}w(\xi_0)|\phi_k(\xi_0)|^2 - 1\right| \le \left|\frac{w(\xi_0)}{\P_n(\xi_0)} - 1\right| + 
\left|\frac{w(\xi_0)}{\F_n(\xi_0)} - 1\right|+64\frac{w(\xi_0)}{\P_n(\xi_0)}\sqrt[4]{\K_n(\xi_0)}.
\end{equation}
\end{Prop}
\beginpf We use inequality \eqref{eq22} from the proof of Theorem \ref{t1},
\begin{equation}\label{eq23}
\frac{1}{n}\sum_{k=0}^{n-1}e^{\P(\log w, z_n)} |\phi_k(\xi)|^2 \le \left(1 - \frac{\delta}{2n}\right)^{-2n} (1 + 2\sqrt{\eps_n} + 4\eps_n),
\end{equation}
where $\delta \in (0,1]$, $z_n = 1 - \delta/n$, and $\eps_n = \frac{e^{\K(\mu, z_n)} -1}{\delta}$. For $x \in  [0, 1/2]$, we have $\log (1-x) \ge -2x$ and $e^{2x} \le 1+4x$, hence 
\begin{align*}
\left(1 - \frac{\delta}{2n}\right)^{-2n} = e^{-2n \log(1- \delta/2n)} \le e^{2 \delta} 
\le 1+ 4 \delta,& \\
2\sqrt{\eps_n} + 4\eps_n \le 6\delta^{-1}\sqrt{e^{\K(\mu, z_n)} - 1} \le 12\delta^{-1}\sqrt{\K_n(\xi_0)}.&
\end{align*}
Choosing $\delta = \sqrt[4]{\K_n(\xi_0)} \le 1$, we  estimate the right hand side in \eqref{eq23} by 
$$
(1 + 4\delta)(1 + 12\delta^{-1}\sqrt{\K_n(\xi_0)}) \le 1 + 64\sqrt[4]{\K_n(\xi_0)}.
$$
This and (Jensen's) inequality $e^{\P(\log w, z_n)} \ge \P_n(\xi_0)$ gives us the upper bound in \eqref{eq26}. The lower bound in \eqref{eq26} is just a combination of estimates \eqref{eq24}, \eqref{eq25}.  \qed

\medskip

Closing this section, let us comment on the continuous version of the problem. Recently P.\,Gubkin \cite{Gub21} proved a variant of Theorem \ref{t1} for Krein's systems. His proof follows the line of A.\,M\'at\'e, P.\,Nevai and V.\,Totik \cite{MNT91}. While the entropy function $\K(\mu, z)$ has sense (and was appeared firstly) in a continuous setting \cite{BD2017}, at the present moment formula \eqref{eq6} has no direct continuous counterpart, preventing an immediate transfer of our arguments to the case of Krein systems or Dirac operators. This remains an interesting open direction.

\medskip

\section{Proof of Theorem \ref{t2}}\label{s3}
Given integers $n_1 \le n_2$, let $L_{n_1, n_2} = \spn\{z^k, \; n_1 \le k \le n_2\}$  be a subspace of trigonometric polynomials on $\T$. For $n \ge 0$, the set $L_{0, n} = P_n$ consists of polynomials of degree at most $n$. In general, elements of $L_{n_1, n_2}$ can be represented as polynomials in two variables, $z$, $\bar z$. Let us denote 
$$
L_n = \spn\{\chi_j, \; 0 \le j \le n\} =
\begin{cases}
L_{-k, k} & n= 2k,\\
L_{-k, k+1} & n = 2k+1.
\end{cases}
$$
We regard $L_n$ as an $n+1$ dimensional subspace of $L^2(\mu)$. For $\xi \in \T$, denote by $k_{L_n, \mu, \xi}$ the reproducing kernel in $L_n$ at $\xi$. By definition, $k_{L_n, \mu, \xi}$ is the element of $L_n$ such that $(f, k_{L_n, \mu, \xi})_{L^2(\mu)} = f(\xi)$ for every $f \in L_n$. It is easy to see that
\begin{equation}\label{eq21}
k_{L_n, \mu, \xi}(z) = \sum_{k=0}^{n}\ov{\chi_k(\xi)}\chi_{k}(z), \qquad z \in \T.
\end{equation}
Our first aim is to derive a variant of Christoffel-Darboux formula for $k_{L_n, \mu, \xi}$. 
\begin{Lem}\label{l5}
Let $\mu$ be a probability measure on $\T$. Let $n \ge 0$ be an integer number, and denote by $k$ the integer part of $n/2$. Then we have
$$
k_{L_n, \mu, \xi}(z) = (\xi\ov{z})^k k_{P_n,\mu,\xi}(z) = 
\begin{cases}
\frac{z\ov{\chi_{n+1}(z)\xi}\chi_{n+1}(\xi) - \chi_{n+1}(z)\ov{\chi_{n+1}(\xi)}}{1 - \ov{\xi}z}, & n \mbox{ is even,}\\
\frac{z\chi_{n+1}(z)\ov{\xi\chi_{n+1}(\xi)} - z\ov{\chi_{n+1}(z)\xi}\chi_{n+1}(\xi)}{1 - \ov{\xi}z}, & n \mbox{ is odd.}
\end{cases}
$$
foe every $\xi,z \in \T$.
\end{Lem}
\beginpf We have $f \in L_n$ if and only if $z^k f \in P_n$. Hence,
$$
\xi^k f(\xi) = (z^k f, k_{P_n, \mu, \xi})_{L^2(\mu)} = ( f, \ov{z}^k k_{P_n, \mu, \xi})_{L^2(\mu)}.
$$
It follows that $k_{L_n, \mu, \xi} = (\xi\ov{z})^k k_{P_n, \mu, \xi}$. 
So, it remains to check the formula for $k_{L_n, \mu, \xi}$ in terms of $\chi_n$. It is well-known that
\begin{equation}\label{eq31}
k_{P_{n}, \mu, \xi} = \frac{\phi_{n+1}^*(z)\ov{\phi_{n+1}^*(\xi)} - \phi_{n+1}(z)\ov{\phi_{n+1}(\xi)}}{1 - \ov{\xi}z},
\end{equation}
see, e.g., Chapter V in Freud \cite{Freudbook}.
Recall that
$$
\chi_{2k} = \bar z^{k}\phi_{2k}^{*}, \qquad  
\chi_{2k+1} = \bar z^{k}\phi_{2k+1}, \qquad z \in \T, \qquad
k \ge 0.
$$
Consider the case $n = 2k$. We have
$$
\ov{z}^k\phi_{n+1}^* = z^{k+1}\ov{\phi_{2k+1}}
=z\ov{\chi_{n+1}}, \qquad \ov{z}^k\phi_{n+1} =  \chi_{n+1}.
$$
Hence,
$$
k_{L_n, \mu, \xi}(z) =  \frac{z\ov{\chi_{n+1}(z)\xi}\chi_{n+1}(\xi) - \chi_{n+1}(z)\ov{\chi_{n+1}(\xi)}}{1 - \ov{\xi}z}
$$
In the case $n = 2k+1$ we have
$$
\ov{z}^k\phi_{n+1}^* = z\chi_{n+1}, \qquad \ov{z}^k\phi_{n+1} = z^{k+2}\ov{\phi_{2k+2}^*} = z\ov{\chi_{n+1}}, $$
therefore,
$$
k_{L_n, \mu, \xi}(z) =  \frac{z\chi_{n+1}(z)\ov{\xi\chi_{n+1}(\xi)} - z\ov{\chi_{n+1}(z)\xi}\chi_{n+1}(\xi)}{1 - \ov{\xi}z}.
$$
The lemma follows. \qed

\medskip

\noindent{\bf Remark.} Another variant of Christoffel-Darboux formula for $k_{L_n,\mu,n}$ was derived by R.\,Cruz-Barroso and P.\,Gonz\'{a}lez-Vera \cite{CG05}.


\medskip

We will use the following well-known fact.

\begin{Lem}\label{l4}
Let $\mu$ be a finite nonnegative measure on $\T$, and let $f \in L^2(\mu)$. Then at every Lebesgue point $\xi_0$ of $f$ we have
$$
\lim_{r\to 1}\frac{\P(|f|^2\,d\mu, r\xi_0)}{\P(\mu, r\xi_0)} = |f(\xi_0)|^2.
$$
\end{Lem}

\medskip

Everything is ready for the proof of Theorem \ref{t2}. As reader will see, the proof essentially uses the original idea of G.\,Freud.

\medskip

\noindent {\bf Proof of Theorem \ref{t2}.}
Take a function $f \in L^2(\mu)$ and assume that $\xi_0 \in \T$ is its Lebesgue point. Denote $I_n = I_{\xi_0, 1/n}$, $E_n = \T\setminus I_n$. Relation \eqref{eq21} gives
\begin{equation}\notag
\bS_n(f, \xi_0) 
= \int_{\T}f\ov{k_{L_n,\mu,\xi_0}}\,d\mu
=\int_{I_n}f\ov{k_{L_n,\mu,\xi_0}}\,d\mu + \int_{E_n}f\ov{k_{L_n,\mu,\xi_0}}\,d\mu = J_{1,n} + J_{2,n}.
\end{equation}
Let us estimate $J_{1,n}$, $J_{2,n}$ separately. By Cauchy-Schwarz inequality, we have
\begin{align*}
\limsup_{n \to +\infty}|J_{1,n}|^2 
&\le \limsup_{n \to +\infty}\int_{I_n}|f|^2\,d\mu \cdot \|k_{L_n,\mu,\xi_0}\|_{L^2(\mu)}^{2} \\
&\le \limsup_{n \to +\infty}\frac{1}{\mu(I_n)}\int_{I_n}|f|^2\,d\mu \cdot \sup_{n \ge 1}\mu(I_n)\|k_{L_n,\mu,\xi_0}\|_{L^2(\mu)}^{2}\\
&\le |f(\xi_0)|^2\cdot \sup_{n \ge 1}\mu(I_n)\|k_{L_n,\mu,\xi_0}\|_{L^2(\mu)}^{2}.
\end{align*}
Set $z_n = (1-1/n)\xi_0$ and note that $\frac{1-|z_n|^2}{|1 - \ov{\xi}z_n|^2} \ge \frac{1 - |z_n|}{(|\xi-\xi_0| + 1/n)^2} \ge n/4$ for $\xi \in I_n$. From this and our assumption we see that
$$
\|k_{L_n,\mu,\xi_0}\|_{L^2(\mu)}^{2} = 
k_{L_n,\mu,\xi_0}(\xi_0) \le \frac{cn}{\P(\mu, z_n)} \le \frac{4c}{\mu(I_n)}.
$$
Therefore, we have
$$
\limsup_{n\to +\infty}\frac{1}{n}\sum_{0}^{n-1}|J_{1,k}|\le  \limsup_{n\to +\infty}|J_{1,n}| \le 2\sqrt{c}|f(\xi_0)|.
$$ 
Turning to $J_{2,n}$, let $\mathbb{I}_{E_n}$ denote the indicator function of $E_n$. Define the function $G = \mathbb{I}_{E_n}f/(1 - \xi_0 \bar\xi)$ and estimate $J_{2,n}$ using Lemma \ref{l5} as follows: 
\begin{align*}
|J_{2,n}| 
\le &|(G, \xi\ov{\chi_{n+1}})_{L^2(\mu)}|\cdot |\chi_{n+1}(\xi_0)| + 
|(G, \chi_{n+1})_{L^2(\mu)}|\cdot |\chi_{n+1}(\xi_0)| + \\  
&+ |(G, \xi\chi_{n+1})_{L^2(\mu)}|\cdot |\chi_{n+1}(\xi_0)|.
\end{align*}
This estimate holds for every $n \ge 1$ regardless its parity. For $\xi \in E_n$, we have
$$
\frac{1}{|1 - \ov{\xi_0}\xi|^2} = \left|\frac{z_n - \xi_0 + \xi_0 - \xi}{ \xi_0 - \xi}\right|^2\frac{1}{|1 - \ov{z_n} \xi|^2} \le 
\frac{4}{|1 - \ov{z_n} \xi|^2} \le 4n\frac{1 - |z_n|^2}{|1 - \ov{z_n}\xi|^2}.
$$
Since $\{\chi_k\}_{k \ge 0}$ is an orthonormal system, we have
\begin{align*}
\left(\sum_{1}^{n}|(G, \chi_k)_{L^2(\mu)}|\cdot |\chi_k(\xi_0)|\right)^2 
&\le \left(\sum_{1}^{n}|(G, \chi_k)_{L^2(\mu)}|^2 \right)\sum_{1}^{n}|\chi_k(\xi_0)|^2,\\
&\le \int_{E_n} \frac{|f(\xi)|^2}{|1-\ov{\xi_0}\xi|^2}\,d\mu(\xi)\cdot \sum_{0}^{n}|\chi_k(\xi_0)|^2,\\
&\le 4n\cdot\P(|f|^2\mu, z_n)\cdot \sum_{0}^{n}|\chi_k(\xi_0)|^2,\\
&\le 4cn^2 \cdot\frac{\P(|f|^2\mu, z_n)}{\P(\mu, z_n)},
\end{align*}
where we used assumption \eqref{eq17}.
From Lemma \ref{l4} we now get
$$
\limsup_{n \to +\infty}\frac{1}{n}\sum_{1}^{n}|(G, \chi_k)_{L^2(\mu)}|\cdot |\chi_k(\xi_0)| \le 2\sqrt{c} \cdot |f(\xi_0)|.
$$
Since $\{\xi\chi_{k}\}_{k \ge 0}$ and $\{\xi\ov{\chi_{k}}\}_{k \ge 0}$ are also orthonormal systems in $L^2(\mu)$, the same reasoning applies and gives us the estimate
$$
\limsup_{n \to +\infty}\frac{1}{n}\sum_{k = 0}^{n-1}|J_{2,n}| \le 6\sqrt{c}\cdot |f(\xi_0)|.
$$
Summarizing, we have proved that
$$
\limsup_{n \to +\infty}\frac{1}{n}\sum_{k = 0}^{n-1}|\bS_n(f, \xi_0)| \le 8\sqrt{c}\cdot |f(\xi_0)|.
$$
Using this bound for $f - f(\xi_0)$ in place of $f$, we obtain 
$$
\limsup_{n \to +\infty}\frac{1}{n}\sum_{k=0}^{n-1}|\bS_n(f, \xi_0) - f(\xi_0)| 
= \limsup_{n \to +\infty}\frac{1}{n}\sum_{k=0}^{n-1}|\bS_n(f-f(\xi_0), \xi_0)| = 0.
$$
In other words, the Fourier series of $f$ is strongly Ces\`aro summable to $f(\xi_0)$ at $\xi_0$ provided $\xi_0$ is a Lebesgue point of $f \in L^2(\mu)$ and assumption \eqref{eq17} holds at $\xi_0$. By Theorem \ref{t1}, this assumption holds Lebesgue almost everywhere on $\T$. This ends the proof. \qed

\medskip

\begin{Prop}\label{p4}
Let $\mu = w\,dm + \mu_s$ be a probability measure on $\T$ such that $w = 0$ or $\log w \in L^1(\T)$, and suppose that $\mu_s$ is purely discrete. Then  
\begin{equation}\label{eq4}
\frac{1}{n}\sum_{k=0}^{n}|\chi_k(\xi_0)|^2 \le \frac{c}{\P(\mu, z_n)}, \qquad z_n = \left(1 - \tfrac{1}{n}\right)\xi_0,
\end{equation}
for $\mu$-almost all $\xi_0 \in \T$ and some constant $c> 0$.
\end{Prop}
\beginpf By Theorem \ref{t1} and the definition of $\chi_k$, \eqref{eq4} holds for $w\,dm$-almost all $\xi \in \T$. So, we can assume that $\xi_0 \in \T$ is such that $\mu_s(\{\xi_0\}) > 0$. Then 
$$
\sum_{k \ge 0}|\chi_k(\xi_0)|^2 < \infty,
$$ 
due to the orthogonality of $\{\chi_k\}_{k \ge 0}$. Therefore, both sides in \eqref{eq4} are comparable to $1/n$ for large $n$, and the result follows\footnote{A similar observation was made at the end of Section 3 in \cite{MN80}.}. \qed

\medskip

Two next propositions are related to Problem \ref{prob2}. 

\begin{Prop}\label{p2}
Let $\mu = w\,dm + \mu_s$ be a probability measure on $\T$, $\log w \in L^1(\T)$, and suppose that $\sum_{n \ge 0}\ov{\phi_n(0)} \phi_n$ converges Lebesgue almost everywhere on $\T$ in Ces\`aro sense. Then for Lebesgue almost every $\xi \in \T$ we have
$$
\lim_{n \to \infty}\frac{1}{n}\sum_{k = 0}^{n-1}|\phi_k^*(\xi)D_{\mu}(\xi) - 1|^2 = 0.
$$	
\end{Prop}
\beginpf In view of Theorem \ref{t1} and relation $|D_\mu|^2 = w$ on $\T$, it suffices to prove that
$$
\lim_{n \to \infty}\frac{1}{n}\sum_{k = 0}^{n-1}\Re(\phi_k^*(\xi)D_{\mu}(\xi)) = 1,
$$
which follows from 
$$
\lim_{n \to \infty}\frac{1}{n}\sum_{k = 0}^{n-1}\phi_k^*(\xi) = D_{\mu}(\xi)^{-1}.
$$
It is known that
\begin{equation}\label{eq10}
\lim_{n \to \infty}\int_{\T}|\phi_k^* - D_\mu^{-1}|^2 w\,dm = 0,
\end{equation}
see Section 2.7 in Simon \cite{Simonbook}.
Therefore, it suffices to prove that for Lebesgue almost every $\xi \in \T$ there exists a limit of functions 
\begin{align*}
\frac{1}{n}\sum_{k = 0}^{n-1}\phi_k^*(\xi)\ov{D_{\mu}^{-1}(0)}
= &\frac{1}{n}\sum_{k = 0}^{n-1}(\phi_k^*(\xi)\overline{\phi_k^*(0)} -  \phi_k(\xi)\overline{\phi_k(0)})  \notag\\
&+\frac{1}{n}\sum_{k = 0}^{n-1}\phi_k(\xi)\overline{\phi_k(0)}
\\
&+ \frac{1}{n}\sum_{k = 0}^{n-1}\phi_k^*(\xi)\overline{(D_{\mu}(0)^{-1}- \phi_k^*(0))}. 
\end{align*}
Since $\mu \in \szc$, we have $\lim_k\phi_k(0) = 0$, $\lim_k \phi_k^*(0)= D_{\mu}(0)^{-1}$, see Section 2.7 in Simon \cite{Simonbook}. By Theorem \ref{t1},  $\sup_{n \ge 0}\frac{1}{n}\sum_{k = 0}^{n-1}|\phi_k(\xi)|^2 < c(\xi) < \infty$ for Lebesgue almost every $\xi \in \T$. Therefore, for Lebesgue almost every $\xi \in \T$ we have
$$
\lim_{n \to \infty}\left|\frac{1}{n}\sum_{k = 0}^{n-1}\phi_k(\xi)\overline{\phi_k(0)}\right|^2 \le \lim_{n \to \infty}\frac{c(\xi)}{n}\sum_{k=0}^{n-1}|\phi_k(0)|^2 = 0,
$$
and
$$
\lim_{n \to \infty}\left|\frac{1}{n}\sum_{k = 0}^{n-1}\phi_k^*(\xi)\overline{(D_{\mu}(0)^{-1}- \phi_k^*(0))} \right|^2 \le \lim_{n \to \infty}\frac{c(\xi)}{n}\sum_{k=0}^{n-1}|D_{\mu}(0)^{-1}- \phi_k^*(0)|^2 = 0,
$$
by regularity of Ces\`aro averaging method. So, we need to prove that Lebesgue almost everywhere on $\T$ there exists the limit
\begin{equation}\label{eq29}
\lim_{n \to \infty}\frac{1}{n}\sum_{k = 0}^{n-1}(\phi_k^*(\xi)\overline{\phi_k^*(0)} -  \phi_k(\xi)\overline{\phi_k(0)}).
\end{equation}
By our assumption, for almost every $\xi \in \T$ there exists the limit
\begin{equation}\label{eq30}
\lim_{n \to \infty}\frac{1}{n}\sum_{k = 0}^{n-1}S_k(\xi), \qquad S_k(\xi) = \sum_{j = 0}^{k}\ov{\phi_j(0)}\phi_j(\xi).
\end{equation}
Observe that $S_k$ is the reproducing kernel in $P_k$ at $0$. By Christoffel-Darboux formula \eqref{eq31}, we have $S_{k-1} = \phi_k^*(\xi)\overline{\phi_k^*(0)} -  \phi_k(\xi)\overline{\phi_k(0)}$. The claim follows. \qed

\medskip

Proposition \ref{p2} has an interesting application in scattering theory. Consider a  probability measure $\mu$ on $\T$. Let $\{\phi_{n}\}_{n \ge 0}$, $\{f_{n}\}_{n \ge 0}$ denote its orthonormal polynomials and Schur family, correspondingly. Put  $a_n = f_n(0)$. It is known (see, e.g, \cite{Simonbook} or \cite{KH01}) that vectors $\Xi_n = 
\left(
\begin{smallmatrix}
\phi_n(z)\\
\phi_n^*(z)
\end{smallmatrix}
\right)$ satisfy the recurrence relation
\begin{equation}\label{eq12}
\Xi_{n+1} = T_{n}\Xi_{n}, \qquad n \ge 0,
\end{equation}
where 
$$
T_n = \rho_n^{-1} \begin{pmatrix}
z & -\ov{a_n}\\
-a_n z & 1
\end{pmatrix}, 
 \qquad \rho_n = \sqrt{1 - |a_n|^2}.  
$$
One approach in the scattering theory of one-dimensional Dirac and Schr\"odinger operators utilizes relations of these operators to their discrete relative -- system \eqref{eq12}, see, e.g., \cite{Den02b}, \cite{B2018}. Following B.\,Simon \cite[Section 10.7]{Simonbook2}, we define Jost solutions, $f_\pm(n,\xi)$, $\xi \in \T$, $n \ge 0$, as the two vectors obeying \eqref{eq12} with
$$
\lim_{n \to \infty}\left\|f_+(n, \xi) - \begin{pmatrix}\xi^n \\ 0\end{pmatrix}\right\| = 0, 
\qquad \lim_{n \to \infty}\left\|f_-(n, \xi) - \begin{pmatrix}0 \\ 1\end{pmatrix}\right\| = 0.
$$
Theorem 10.7.7 in B.\,Simon \cite{Simonbook2} gives the existence of Jost solutions for every $\xi \in \T$ under the Baxter's condition $\sum_{n \ge 0}|a_n| < \infty$. This corresponds to the very classical case of potentials in $L^1(\R)$ for Dirac and Schr\"odinger operators. Much more involved case $V \in L^p(\R)$, $1 \le p < 2$, treated by M.\,Christ and A.\,Kiselev in \cite{ChK01} corresponds to condition $\sum_{n \ge 0}|a_n|^p < \infty$. Existence of Jost functions for \eqref{eq12} for almost every $\xi \in \T$ in this case is still an open problem. Let us show that positive answer on Problem \ref{prob2} implies the existence of averaged Jost functions for every sequence $\{a_n\}$ such that $\sum_{n \ge 0}|a_n|^2 < \infty$. For this we need the notion of a dual orthogonality measure. 

\medskip

Let $\mu$ be a probability measure on $\T$, and let $f$ be its Schur function. The probability measure $\nu$ whose Schur function equals $-f$ is called the dual orthogonality measure. It can be defined by the relation
$$
\frac{1 - |zf(z)|^2}{|1+ zf(z)|^2} = \P(\nu, z), \qquad z \in \D. 
$$ 
We have $\mu \in \szc$ if and only if $\nu \in \szc$ if and only if $\log(1 - |f|^2) \in L^1(\T)$. Therefore, the Szeg\H{o} function $D_{\nu}$ is defined. We also will use notation $\psi_n$, $\psi_n^*$ for the orthonormal polynomials with respect to $\nu$ and their reflected versions, $\psi_n^* = z^n \ov{\psi_n(1/\bar z)}$. 
\begin{Prop}\label{p3}
Let $\mu \in \szc$, and let $\nu$ be the dual orthogonality measure with respect to $\mu$. Suppose that series $\sum_{n \ge 0}\ov{\phi_n(0)} \phi_n$, $\sum_{n \ge 0}\ov{\psi_n(0)} \psi_n$ converge Lebesgue almost everywhere on $\T$ in Ces\`aro sense. Then for Lebesgue almost every $\xi \in \T$ there exist solutions $f_{\pm}(n,\xi)$ of \eqref{eq12} such that 
\begin{equation}\label{eq11}
\lim_{n \to \infty}\frac{1}{n}\sum_{k=0}^{n-1}\left\|f_+(k, \xi) - \begin{pmatrix}\xi^n \\ 0\end{pmatrix}\right\| = 0, 
\qquad \lim_{n \to \infty}\frac{1}{n}\sum_{k=0}^{n-1}\left\|f_-(k, \xi) - \begin{pmatrix}0 \\ 1\end{pmatrix}\right\| = 0.
\end{equation}
\end{Prop}
\beginpf Following B.\,Simon (see Section 10.7 in \cite{Simonbook2}), we define $F = D_{\mu}/D_{\nu}$ and
$$
f_+(n, \xi) = \frac{1}{2}D_{\mu}^{-1}(\xi)\left[\begin{pmatrix}\psi_n(\xi)\\ -\psi_n^*(\xi)\end{pmatrix} + F(\xi)\begin{pmatrix}\phi_n(\xi)\\ \phi_n^*(\xi)\end{pmatrix}\right].
$$
It is clear that $f_+(n,\xi)$ is a solution of \eqref{eq12}. 
Almost everywhere $\T$ we have
$$
\begin{pmatrix}\phi_n\\ \phi_n^*\end{pmatrix} 
= 
\begin{pmatrix}\xi^n\ov{D_{\mu}^{-1}}\\ D_{\mu}^{-1}\end{pmatrix} 
+
\begin{pmatrix}\xi^n(\ov{\phi_n^* - D_{\mu}^{-1}})\\ \phi_n^* - D_{\mu}^{-1}\end{pmatrix}, 
$$
and
$$
\begin{pmatrix}\psi_n\\ -\psi_n^*\end{pmatrix} 
= 
\begin{pmatrix}\xi^n\ov{D_{\nu}^{-1}}\\ -D_{\nu}^{-1}\end{pmatrix} 
+
\begin{pmatrix}\xi^n(\ov{\psi_n^* - D_{\nu}^{-1}})\\ -\psi_n^* + D_{\nu}^{-1}\end{pmatrix}.
$$
By Proposition \ref{p2}, for almost every $\xi \in \T$ there exist sequences $\eps_{1,n}(\xi)$, $\eps_{2,n}(\xi)$ converging to zero in the Ces\`aro sense and such that 
$$
f_+(n,\xi) = \frac{1}{2}D_{\mu}^{-1}(\xi)
\left[\begin{pmatrix}\xi^n(\ov{D_{\nu}^{-1}(\xi)} + F(\xi)\ov{D_{\mu}^{-1}(\xi)}) \\
-D_{\nu}^{-1}(\xi) + F(\xi)D_{\mu}^{-1}(\xi) 
\end{pmatrix} + \begin{pmatrix}\eps_{1,n}(\xi)\\ \eps_{2,n}(\xi)\end{pmatrix}\right].
$$
Since $F = D_{\mu}/D_{\nu}$, we have
$$
f_+(n,\xi) = \begin{pmatrix}\xi^n\bigl(\Re (D_{\mu}^{-1}(\xi)\ov{D_{\nu}^{-1}(\xi)}\bigr) \\
0
\end{pmatrix} + \begin{pmatrix}\eps_{1,n}(\xi)\\ \eps_{2,n}(\xi)\end{pmatrix}.
$$
Let us show that $\Re \bigl(D_{\mu}^{-1}\ov{D_{\nu}^{-1}}\bigr) = 1$ almost everywhere on $\T$ (then the first relation in \eqref{eq11} will follow). We have $\Re\bigl(\phi_n^*(\xi)\ov{\psi_n^*(\xi)}\bigr) = 1$ for every $n \ge 0$ and $\xi \in \T$, see formula (3.2.22) in \cite{Simonbook}. Note that $\phi_n^* \to D_\mu^{-1}$,  $\psi_n^* \to D_\nu^{-1}$ in Lebesgue measure on $\T$ by \eqref{eq10}. Choosing a pointwise convergent subsequence on a set of full Lebesgue measure, we obtain identity $\Re \bigl(D_{\mu}^{-1}\ov{D_{\nu}^{-1}}\bigr) = 1$  and the first relation in \eqref{eq11}. Similarly, one can check that 
$$
f_-(n, \xi) = -\frac{1}{2}\ov{D_{\mu}^{-1}(\xi)}\left[\begin{pmatrix}\psi_n(\xi)\\ -\psi_n^*(\xi)\end{pmatrix} - \ov{F(\xi)}\begin{pmatrix}\phi_n(\xi)\\ \phi_n^*(\xi)\end{pmatrix}\right].
$$  
is a solution of \eqref{eq12} satisfying the second relation in \eqref{eq11}. Indeed, we have
$$
f_-(n, \xi) = 
-\frac{1}{2}\ov{D_{\mu}^{-1}(\xi)}\left[\begin{pmatrix}\xi^n(\ov{D_{\nu}^{-1}(\xi) - F(\xi) D_{\mu}^{-1}(\xi)})\\ -D_{\nu}^{-1}(\xi) - \ov{F(\xi)} D_{\mu}^{-1}(\xi)\end{pmatrix}  + \begin{pmatrix}\tilde\eps_{1,n}(\xi)\\ \tilde\eps_{2,n}(\xi)\end{pmatrix} \right],
$$
for some sequences $\tilde\eps_{1,n}(\xi)$, $\tilde\eps_{2,n}(\xi)$ converging to zero in the Ces\`aro sense Lebesgue  almost everywhere on $\T$.
Since $D_{\nu}^{-1} - FD_{\mu}^{-1} = 0$ and
$$
D_{\nu}^{-1} + \ov{F} D_{\mu}^{-1} = 
2D_{\mu}^{-1}\Re(F) = 2D_{\mu}^{-1}|D_\mu|^{2} = 2\ov{D_\mu},  
$$
almost everywhere on $\T$, the result follows. \qed

\medskip

As a final remark, let us mention that recently   A.\,Poltoratskii \cite{Polt21} established the pointwise convergence almost everywhere on the real line of absolute values of solutions of Dirac systems with $L^2(\R)$--potentials. In the theory of orthogonal polynomials his result ``corresponds'' to the following (still unproved) assertion: for every $\mu \in \szc$, we have $\lim_{n \to \infty}|\phi_n(\xi)| \to |D_{\mu}^{-1}(\xi)|$ Lebesgue almost everywhere on $\T$.



\bibliographystyle{plain} 
\bibliography{bibfile}
\enddocument